\documentclass[a4paper]{amsart}

\usepackage{graphicx,color,tikz,comment}

\theoremstyle{plain}
\newtheorem{thm}{Theorem}

\theoremstyle{definition}

\theoremstyle{remark}

\def\h#1#2#3{\draw[black,line width=0.75pt, cap=round] (#1,#2)--(#1+#3,#2);}
\def\w{0.2}
\def\v#1#2#3{
       \draw[white,line width=2.8pt] (#1,#2+\w)--(#1,#2+#3-\w);
       \draw[black,line width=0.75pt, cap=round] (#1,#2)--(#1,#2+#3);}
\def\vv#1#2#3{
       \draw[black,line width=0.75pt, cap=round] (#1,#2)--(#1,#2+#3);}
\def\hr#1#2#3{\draw[red,line width=0.95pt, cap=round] (#1,#2)--(#1+#3,#2);}
\def\vr#1#2#3{
       \draw[white,line width=2.5pt] (#1,#2+\w)--(#1,#2+#3-\w);
       \draw[red,line width=0.95pt, cap=round] (#1,#2)--(#1,#2+#3);}
\def\hw#1#2#3{\draw[white,line width=0.95pt, cap=round] (#1,#2)--(#1+#3,#2);}

\title[Minimal grid diagrams of  14 crossing knots with arc index 13, 14]{
Minimal grid diagrams of the prime knots with crossing number 14 and arc index 13, 14}

\author[H. J. Lee et al.]{Hwa Jeong Lee}
\address{Dongguk University WISE}
\email{hjwith@dongguk.ac.kr}

\author[]{Alexander Stoimenow}
\address{Dongguk University WISE}
\email{stoimeno@stoimenov.net}

\author[]{Hun Kim}
\address{Korea Science Academy of KAIST}
\email{hunkim@ksa.kaist.ac.kr}

\author[]{Minchae Kim$^\dagger$}
\address{Korea Science Academy of KAIST}
\email{22-017@ksa.hs.kr}

\author[]{Songwon Ryu$^\dagger$}
\address{Korea Science Academy of KAIST}
\email{songwon0425@gmail.com}

\author[]{Dongju Shin$^\dagger$}
\address{Korea Science Academy of KAIST}
\email{djshin7233@gmail.com}

\author[]{Joon Hyeok Choi$^\star$}
\address{Korea Science Academy of KAIST}
\email{\
23-121@ksa.hs.kr}
\author[]{Woo Jin Choi$^\star$}
\address{Korea Science Academy of KAIST}
\email{
23-118@ksa.hs.kr}
\author[]{Jin Seong Park$^\star$}
\address{Korea Science Academy of KAIST}
\email{
23-057@ksa.hs.kr}

\author[]{Gyo Taek Jin}
\address{KAIST}
\email{trefoil@kaist.ac.kr}

\date{\today}
\makeatletter
\@namedef{subjclassname@2020}{\textup{2020} Mathematics Subject Classification}
\makeatother

\begin{document}

\begin{abstract}
There are 46,972 prime knots with crossing number 14. Among them 19,536 are
alternating and have arc index 16. Among the non-alternating knots, 17,
477, and 3,180 have arc index 10, 11, and 12, respectively. The remaining
23,762 have arc index 13 or 14. There are none with arc index smaller than 10 or larger than 14.
 We obtained 8,027 knots having arc index 13 and 15,735 knots having arc index 14. We
show them by their minimal grid diagrams. 
\end{abstract}

\maketitle

\section{Introduction}
\newcommand{\mythicklines}{\linethickness{1.0pt}}
\newcommand{\ArcPreGrid}{
\thinlines\color{cyan}
\put(45,20){\vector(0,1){190}}
\multiput(0,175)(5,0){24}{\line(1,0){2}} 
\multiput(60,145)(5,0){18}{\line(1,0){2}}
\multiput(30,115)(5,0){12}{\line(1,0){2}}
\multiput(0,85)(5,0){12}{\line(1,0){2}}
\multiput(30,55)(5,0){18}{\line(1,0){2}}
\multiput(90,25)(5,0){12}{\line(1,0){2}}
\color{black}\thinlines
\qbezier(0,175)(45,190)(45,190)\qbezier(45,190)(45,190)(120,175)
\qbezier(60,145)(45,160)(45,160)%
\qbezier(45,160)(45,160)(115,150)\qbezier(125,148)(125,148)(150,145)
\qbezier(30,115)(45,130)(45,130)%
\qbezier(45,130)(45,130)(55,127)\qbezier(65,123)(65,123)(90,115)
\qbezier(0,85)(25,93)(25,93)\qbezier(35,97)(45,100)(45,100)%
\qbezier(45,100)(45,100)(60,85)
\qbezier(30,55)(45,70)(45,70)%
\qbezier(45,70)(45,70)(85,62)\qbezier(95,60)(95,60)(120,55)
\qbezier(90,25)(45,40)(45,40)%
\qbezier(45,40)(45,40)(85,34)\qbezier(95,33)(95,33)(150,25)
}

\setlength{\unitlength}{0.015cm}

\begin{figure}[h!]%
\centerline{%
\begin{picture}(80,250)(0,-30)
\put(0,0){\vector(0,1){210}} %
\put(10,-20){\scriptsize$\theta=0$}
\thicklines
\put(0,85){\line(1,0){70}}
\put(70,85){\line(0,1){90}}
\put(0,175){\line(1,0){70}}
\end{picture}\quad
\begin{picture}(80,250)(0,-30)
\put(0,0){\vector(0,1){210}} %
\put(10,-20){\scriptsize$\theta=\frac{\pi}{10}$}
\thicklines
\put(0,55){\line(1,0){70}}
\put(70,55){\line(0,1){60}}
\put(0,115){\line(1,0){70}}
\end{picture}\quad
\begin{picture}(80,250)(0,-30)
\put(0,0){\vector(0,1){210}} %
\put(10,-20){\scriptsize$\theta=\frac{\pi}{5}$}
\thicklines
\put(0,85){\line(1,0){70}}
\put(70,85){\line(0,1){60}}
\put(0,145){\line(1,0){70}}
\end{picture}\quad
\begin{picture}(80,250)(0,-30)
\put(0,0){\vector(0,1){210}} %
\put(10,-20){\scriptsize$\theta=\frac{3\pi}{10}$}
\thicklines
\put(0,25){\line(1,0){70}}
\put(70,25){\line(0,1){90}}
\put(0,115){\line(1,0){70}}
\end{picture}\quad
\begin{picture}(80,250)(0,-30)
\put(0,0){\vector(0,1){210}} %
\put(10,-20){\scriptsize$\theta=\frac{2\pi}{5}$}
\thicklines
\put(0,55){\line(1,0){70}}
\put(70,55){\line(0,1){120}}
\put(0,175){\line(1,0){70}}
\end{picture}\quad
\begin{picture}(80,250)(0,-30)
\put(0,0){\vector(0,1){210}} %
\put(10,-20){\scriptsize$\theta=\frac{\pi}{2}$}
\thicklines
\put(0,25){\line(1,0){70}}
\put(70,25){\line(0,1){120}}
\put(0,145){\line(1,0){70}}
\end{picture}\qquad
\begin{picture}(150,250)(0,-30)
\thicklines
\put(0,85){\line(0,1){90}}
\put(30,55){\line(0,1){60}}
\put(60,85){\line(0,1){60}}
\put(90,25){\line(0,1){90}}
\put(120,55){\line(0,1){120}}
\put(150,25){\line(0,1){120}}
\ArcPreGrid
\end{picture}
}
\caption{An arc presentation of the figure eight knot}
\end{figure}
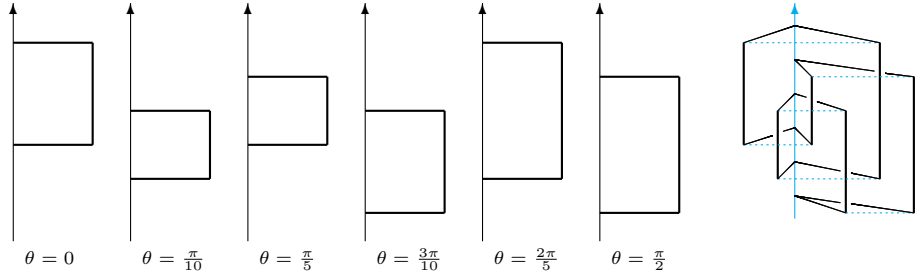

An \emph{arc presentation\/} of a knot $K$ in $\mathbb R^3$ is an embedding of $K$ on a  finite union of half planes which have a common boundary axis such that each half plane contains a properly embedded single arc of $K$.  Such half planes are called \emph{pages}. The minimal number of pages among all arc presentations of $K$ is called the \emph{arc index\/} of $K$, 
denoted by $\alpha(K)$~\cite{C1995}.

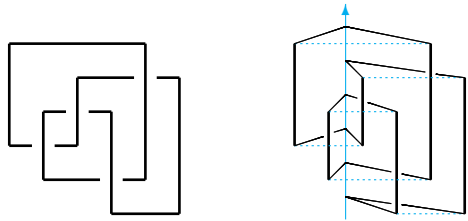
\begin{figure}[h]
\begin{picture}(150,180)(0,10)
\mythicklines
\put(0,85){\line(0,1){90}}
\put(30,55){\line(0,1){60}}
\put(60,85){\line(0,1){60}}
\put(90,25){\line(0,1){90}}
\put(120,55){\line(0,1){120}}
\put(150,25){\line(0,1){120}}
\put(0,175){\line(1,0){120}}
\put(60,145){\line(1,0){50}} \put(150,145){\line(-1,0){20}}
\put(30,115){\line(1,0){20}} \put(90,115){\line(-1,0){20}}
\put(0,85){\line(1,0){20}} \put(60,85){\line(-1,0){20}}
\put(30,55){\line(1,0){50}} \put(120,55){\line(-1,0){20}}
\put(90,25){\line(1,0){60}}
\end{picture}
\qquad\qquad%
\begin{picture}(150,180)(0,10)
\mythicklines
\put(0,85){\line(0,1){90}}
\put(30,55){\line(0,1){60}}
\put(60,85){\line(0,1){60}}
\put(90,25){\line(0,1){90}}
\put(120,55){\line(0,1){120}}
\put(150,25){\line(0,1){120}}
\ArcPreGrid
\end{picture}
\caption{A grid diagram and its corresponding arc presentation}
\label{grid2arcs}
\end{figure}

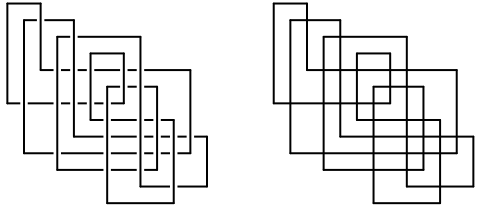
\begin{figure}[h!]
\begin{tikzpicture}[scale=0.2200]
\h{1}{13}{2}
\h{3}{9}{9}
\h{2}{4}{10}
\h{2}{12}{3}
\h{5}{5}{8}
\h{9}{2}{4}
\h{4}{11}{5}
\h{4}{3}{6}
\h{7}{8}{3} 
\h{7}{1}{4}
\h{6}{6}{5}
\h{6}{10}{2}
\h{1}{7}{7}
\v{1}{7}{6}
\v{3}{9}{4}
\v{12}{4}{5}
\v{2}{4}{8} 
\v{5}{5}{7} 
\v{13}{2}{3}
\v{9}{2}{9} 
\v{4}{3}{8} 
\v{10}{3}{5}
\v{7}{1}{7} 
\v{11}{1}{5}
\v{6}{6}{4} 
\v{8}{7}{3} 
\end{tikzpicture}
\qquad
\begin{tikzpicture}[scale=0.2200]
\h{1}{13}{2}
\h{3}{9}{9}
\h{2}{4}{10}
\h{2}{12}{3}
\h{5}{5}{8}
\h{9}{2}{4}
\h{4}{11}{5}
\h{4}{3}{6}
\h{7}{8}{3} 
\h{7}{1}{4}
\h{6}{6}{5}
\h{6}{10}{2}
\h{1}{7}{7}
\vv{1}{7}{6}
\vv{3}{9}{4}
\vv{12}{4}{5}
\vv{2}{4}{8} 
\vv{5}{5}{7} 
\vv{13}{2}{3}
\vv{9}{2}{9} 
\vv{4}{3}{8} 
\vv{10}{3}{5}
\vv{7}{1}{7} 
\vv{11}{1}{5}
\vv{6}{6}{4} 
\vv{8}{7}{3} 
\end{tikzpicture}
\caption{Minimal grid diagrams of the knot $14n10$}\label{fig:14n10tikz}
\end{figure}

A \emph{grid diagram} is a knot diagram with finitely many vertical segments and the same number of  horizontal segments such that the vertical segments cross over the horizontal segments at all crossings. 
Every knot can be drawn as a grid diagram. 
Pulling each horizontal segment in the middle so that it touches an axis behind, we get an arc presentation as shown in Figure~\ref{grid2arcs}. The arc index can be defined as the minimal number of vertical segments (equivalently horizontal segments) among all grid diagrams of a given knot.
As shown in Figure\,\ref{fig:14n10tikz} we sometimes draw grid diagrams without breaking underpassing arcs.

If $D$ is a diagram of a knot or a link $L$,  then $c(D)$ and $c(L)$ denote the number of crosssings in $D$ and the minimal number of crossings in all diagrams of $L$, respectively.

\begin{thm}[\cite{Jin-Park2010}]\label{not bigger than xing}
A prime link $L$ is non-alternating if and only if
$$\alpha(L)\le c(L).$$
\end{thm}

\begin{thm}[\cite{BP2000}]\label{BP2000}
For a non-split link $L$, $\alpha(L)\le c(L)+2$.
\end{thm}

\begin{thm}[\cite{MB1998}]\label{MB1998}
For every link $L$, the arc index  $\alpha(L)$ is bigger than or equal to the $a$-deree span of the Kauffman polynomial $F(a,z)$ plus 2.
\end{thm}

\begin{table}[b]
\small
\caption{Prime knots up to crossing number 14}\label{table1}
{
\newcommand{\vsp}{\vrule width0pt height 10pt depth.5pt} 
\begin{tabular}{|c|c|c|c|c|c|c|c|c|c|c|c|c|c|}
\hline
\vrule width0pt height 12pt depth22pt
\lower15pt\hbox{Crossings} \kern-35pt Arc index &
        \hbox to 7pt{\hfil 5\hfil}&
        \hbox to 7pt{\hfil 6\hfil}&
        \hbox to 7pt{\hfil 7\hfil}&
        \hbox to 7pt{\hfil 8\hfil}&
        \hbox to 7pt{\hfil 9\hfil}&
        \hbox to 7pt{\hfil {10}\hfil}&
        \hbox to 7pt{\hfil {11}\hfil}&
        \hbox to 7pt{\hfil {12}\hfil}&
        \hbox to 7pt{\hfil {13}\hfil}&
        \hbox to 7pt{\hfil {14}\hfil}&
        \hbox to 7pt{\hfil {15}\hfil}& 
        \hbox to 7pt{\hfil {16}\hfil}& Total\\
\hline
\vsp  3&1 & & & & &  &  &  &  &  &  &  &  1\\
\hline
\vsp  4& &1 & & & &  &  &  &  &  &  &  &  1\\
\hline
\vsp  5& & &2 & & &  &  &  &  &  &  &  &  2\\
\hline
\vsp  6& & & &3 & &  &  &  &  &  & &   &  3\\
\hline
\vsp  7& & & & &7 &  &  &  &  &  &  &  &  7\\
\hline
\vsp  8& & &1&2& &{18} &  &  &  & &   &  &21\\
\hline
\vsp  9& & & &2&6&  &{41} &  &  &  & &   & 49\\
\hline
\vsp 10& & & &1&9&32&  &  123 &  & &   &  &  165\\
\hline
\vsp 11& & & & &4&46&135&  & 367 &  &  &  & 552\\
\hline
\vsp 12& & & & &2&48&  211&   627&  &1288 &  &  & 2176\\
\hline
\vsp 13& & & & & &49&  399 & 1412 &  3250 & &  4878&   &9988\\
\hline
\vsp 14& & & & & &17&  477  & 3180 &\bf 8027 & \bf 15735&& 19536 &46972\\
\hline
\end{tabular}\\
}
\end{table}

Table~\ref{table1} shows the number of prime knots of given crossing number and arc index, up to 14 crossings~\cite{HTW1998}. 
Minimal grid diagrams of knots and links up to arc index 9 first appear in~\cite{N1999}. 
The authors of \cite{Jin2006}, \cite{Jin-Park2010,Jin-Park2011}, and \cite{Jin2020,Jin2021} generated Cromwell matrices to find minimal grid diagrams of prime knots with arc index 10, 11, and 12, respectively.

Theorems~\ref{BP2000} and \ref{MB1998} together show that the diagonal entries of Table~\ref{table1} count alternating knots~\cite{HTW1998, Jin-Lee2020, Jin-Lee2022}.  
Theorem~\ref{not bigger than xing} explains why there are none one below the main diagonal entries.

There are 46,972 prime knots with crossing number 14 including 19,536  alternating ones~\cite{HTW1998}. The  alternating ones have arc index 16. Among the other non-alternating knots, 17, 477, and 3,180  have arc index 10, 11,  and 12, respectively. The remaining 23,762 knots have arc index 13 or 14. There are none with arc index smaller than 10 or bigger than 14.
Among the 23,762 knots, 15,735 have the $a$-degree span of the Kauffman polynomial $F_K(a,z)$ equal to 12~\cite{knotscape}. By Theorem~\ref{MB1998}, these 15,735 knots have arc index 14, and thus those having arc index 13 are among the 8,027 remaining.
In this work, we obtained grid diagrams  of size 13 for the 8027 knots and thus showing that their arc index is 13.

\section{The knot-spoke method and the filtered spanning tree method}
The method of filtered spanning tree is a procedure based on the method of knot-spoke diagrams~\cite{BP2000,Jin-Lee2012}. 
A \emph{filtered spanning tree\/} on a prime knot diagram is a sequence of edges\footnote{An edge of a diagram is an arc between two neighboring crossings} forming a maximal tree such that
\begin{itemize}
\item each partial sequence is connected, i.e., is a tree,
\item each partial sequence must not separate untouched crossings into two parts,
\item except at the last step avoid choosing an edge whose one-step extension in the same direction makes a loop. 
\end{itemize}

The proof of Theorem~\ref{BP2000} shows that a filtered spanning tree exists on a prime\footnote{Not a connected sum} diagram.
Let $N$ be a tubular neighborhood of a filtered spanning tree of a diagram $D$ enclosing all crossings inside. Then there are $c(D)+1$ disjoint arcs properly embedded outside $N$.  We deform the arc which is the extension of the last edge of the filtered spanning tree so that a midpoint  touches the boundary of $N$ without intersecting the other arcs.  Inside $N$,  each arc is given a constant height so that the  arcs outside $N$ can be deformed to be placed in vertical half planes to form an arc presentation when the cylinder over $N$ is shrunk to the common axis of the half planes. In this procedure, we obtain an arc presentation in $c(D)+2$ pages. We need to reduce  the number by 3 to get arc index 13 from 14 crossing knots. We must choose the filtered spanning trees so that the number of pages can be reduced by 3. We show this by an example.

\section{A construction using $14n10$}
Every non-alternating prime diagram of a  knot   has at least 4 non-alternating edges which are located on the boundaries of at least 4 regions. 
In the proof of Theorem~\ref{not bigger than xing} in~\cite{Jin-Park2010},
two of the non-alternating edges are used to reduce the number of pages by 2.

\begin{figure}[t]
\centerline{%
\includegraphics[height=0.4\textwidth]{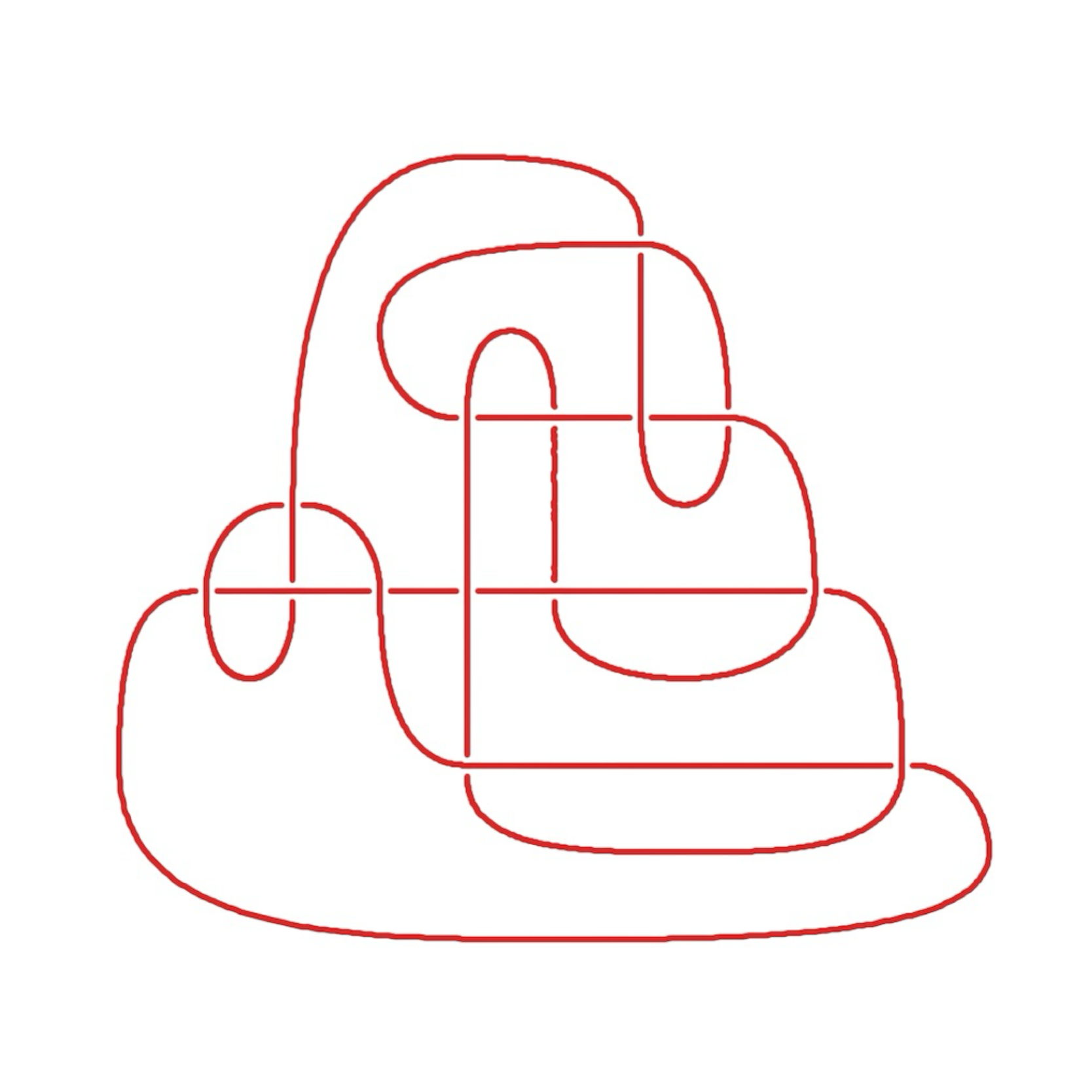}\quad
\includegraphics[height=0.4\textwidth]{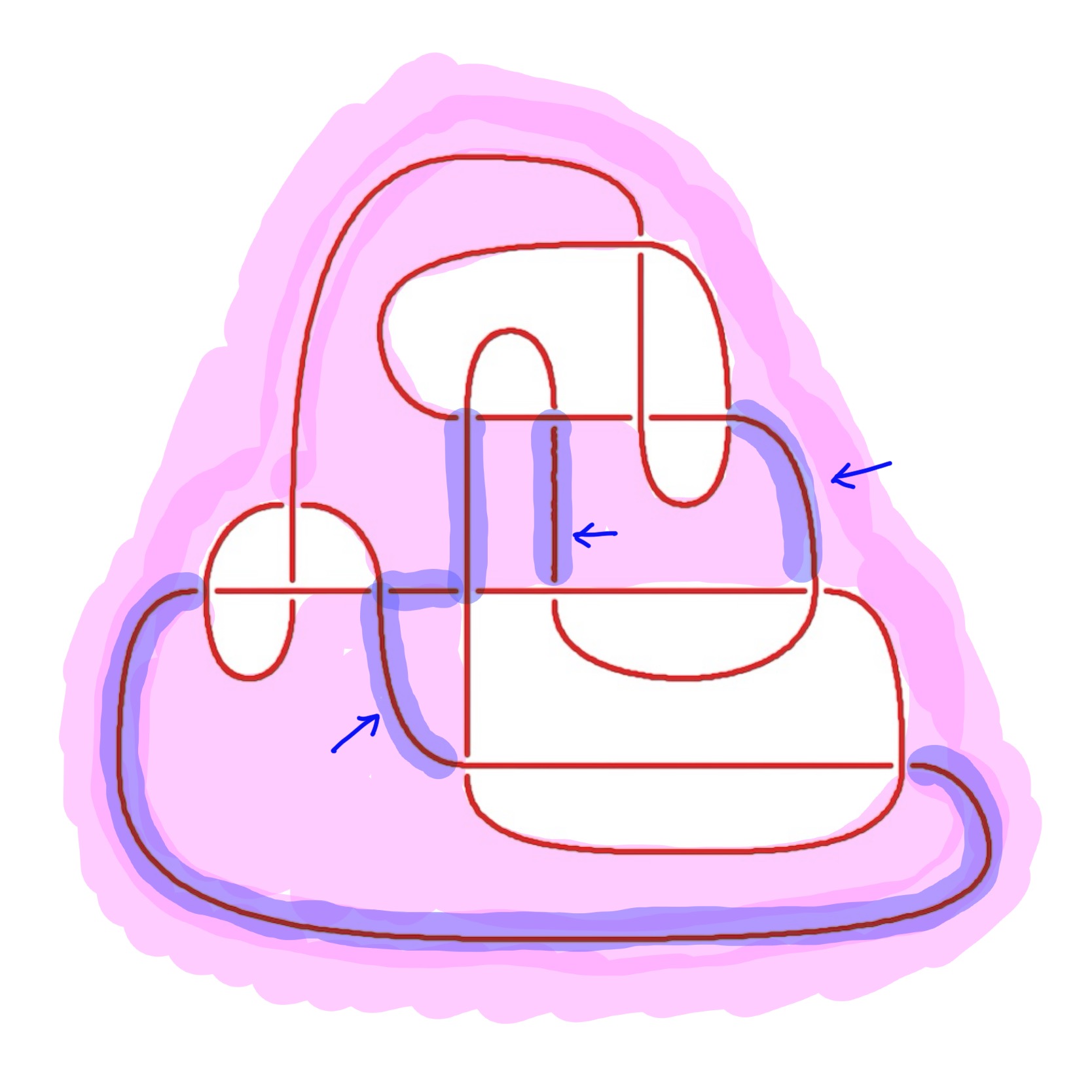}}
\caption{Diagram of $14n10$ and its non-alternating edges \cite{SnapPy}}\label{fig:14n10edges}
\end{figure}

The diagram of $14n10$ in Figure~\ref{fig:14n10edges} has 6 non-alternating edges. We use the three arrowed non-alternating edges to get arc index 13 from the arc presentation of 16 pages obtained from the method of filtered spanning tree.

\begin{figure}[t]
\centerline{%
\includegraphics[height=0.4\textwidth]{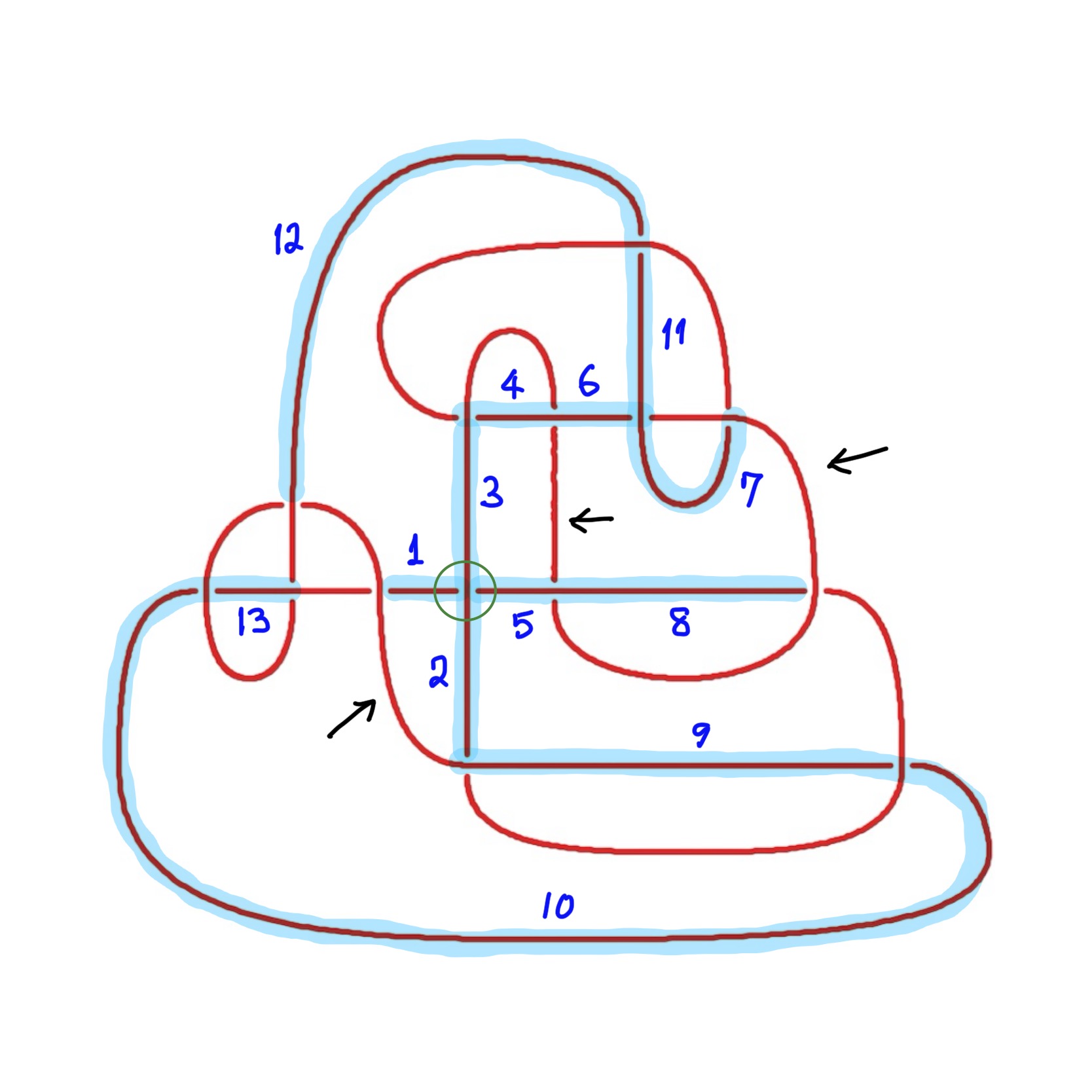}
\includegraphics[height=0.4\textwidth]{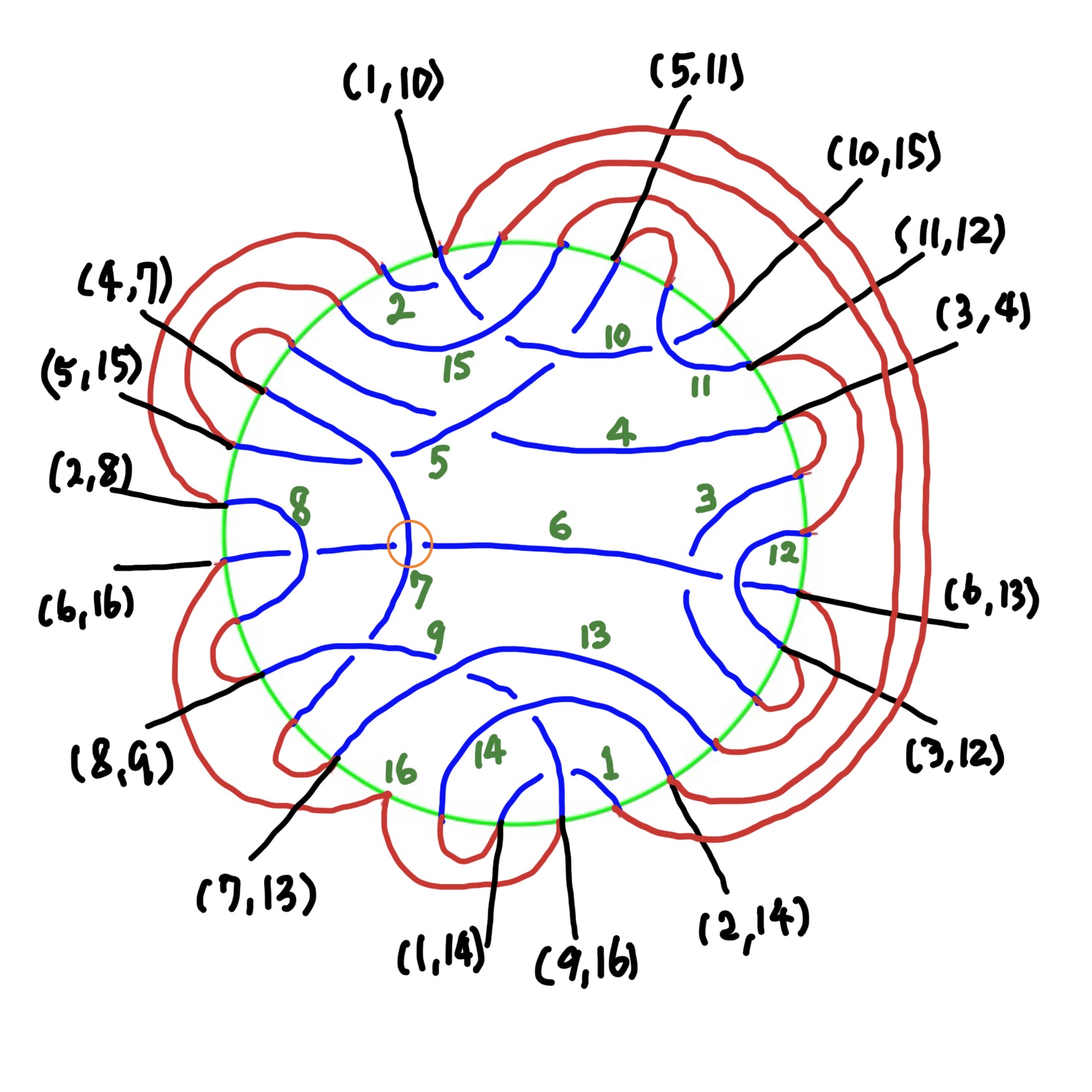}}
\caption{A filtered spanning tree and the resulting spokes}
\label{tree-and-spokes}
\end{figure}

The left of Figure~\ref{tree-and-spokes} shows a filtered spanning tree on $14n10$ where the numbers indicate the order of the edges of the tree. 
The right part is isotoped from the left part so that $N$ becomes the circular disk.  Inside $N$ there are 15 arcs labeled 1 through 15.  These numbers indicate the vertical heights of the arcs in the cylinder over $N$.    
The heights of the arcs inside $N$ are given by the following steps:
\begin{itemize}
\item The arc of the first edge is given an arbitrary height.
\item As one adds a new edge for the filtered spanning tree, the new height is one level higer than the maximun height or one level lower than the minimum height depending on the crossing where the new edge is attached. If it is an extension of an existing edge, no new height is given.
\item If adding an edge $e$ not in the filtered spanning tree makes a cycle, then we give two different heights to the ends of $e$ in a similar way as the above step.
\end{itemize}

Assigning heights as above determines the heights of the end points of the arcs outside $N$ except the midpoint of the extension of the last edge of the filtered spanning tree. We give it one  higer than the maximun height or one  lower than the minimum height. Finally we shift the heights so that they become the numbers 1 through 16.
\bigskip

For each innermost arc outside $N$, we place a spoke\footnote{At the end of the construction, spokes are the projections of the pages onto a plane perpendicular to the axis.} labeled with two heights of its end points, say $(a,b)$, at any one of its two end points. In the figure, they are
$$(1,14), (7,13), (8,9), (4,7), (5,11), (3,4), (3,12).
$$

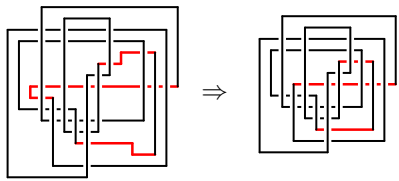
\begin{figure}[h]
\begin{tikzpicture}[scale=0.15]
\input tikz14n10a.tex
\end{tikzpicture}
\begin{picture}(40,150)(10,10)
\put(20,80){$\Rightarrow$}
\end{picture}
\begin{tikzpicture}[scale=0.15]
 \draw[white,line width=2.8pt] (0,-1)--(15,15);
\input tikz14n10b.tex
\end{tikzpicture}

\caption{Destabilization at three places}\label{16grid to 13grid}
\end{figure}

For each arc enclosing other arcs, if one end, say $x$, has its height between two different heights of the inner arcs and the other doesn't, then we place a spoke at~$x$. If both ends are not between  two heights of any inner arcs, we place a spoke at any ends. Such placement of spokes guarantees that after shrinking the cylinder over $N$ the outside arcs can be pushed into the half pages lying over the spokes.
Reading the spokes of Figure~\ref{tree-and-spokes} clockwise starting from $(1,14)$, we have the sequence of vertical intervals

$$\begin{aligned}
&[1,14], [7,13], [8,9],[6,16],[2,8],[5,15],[4,7],[1,10], \\
&[5,11],[10,15],[11,12],[3,4],[6,13],[3,12],[2,14],[9,16].
\end{aligned}$$
In the $xy$-plane, we draw vertical segments $[1,14], [7,13], [8,9], \ldots$, on the vertical  lines $x=1,2,3,\ldots$, respectively. We add horizontal segments joining end points of the vertical segments to obtain the grid diagram on the left of 
Figure~\ref{16grid to 13grid}.

The three arrowed non-alternating edges in Figure~\ref{fig:14n10edges}
correspond to the spokes $(8,9)$, $(11,12)$, and $(3,4)$. As Figure~\ref{16grid to 13grid} illustrates, we can remove their corresponding vertical segments,  and thus obtain a grid diagram of size 13.

\section{Four stages of the work}
Among the 8,027 Dowker-Thistlethwaite codes obtained from knotscape, only 6,295 of them have diagrams with 6 or more non-alternating edges.
Using the method described above, we've sorted out 3,353 knots  having arc index 13. This was the first stage.

\begin{figure}[h]
\begin{picture}(325,140)
\put(0,0){\includegraphics[width=0.4\textwidth]{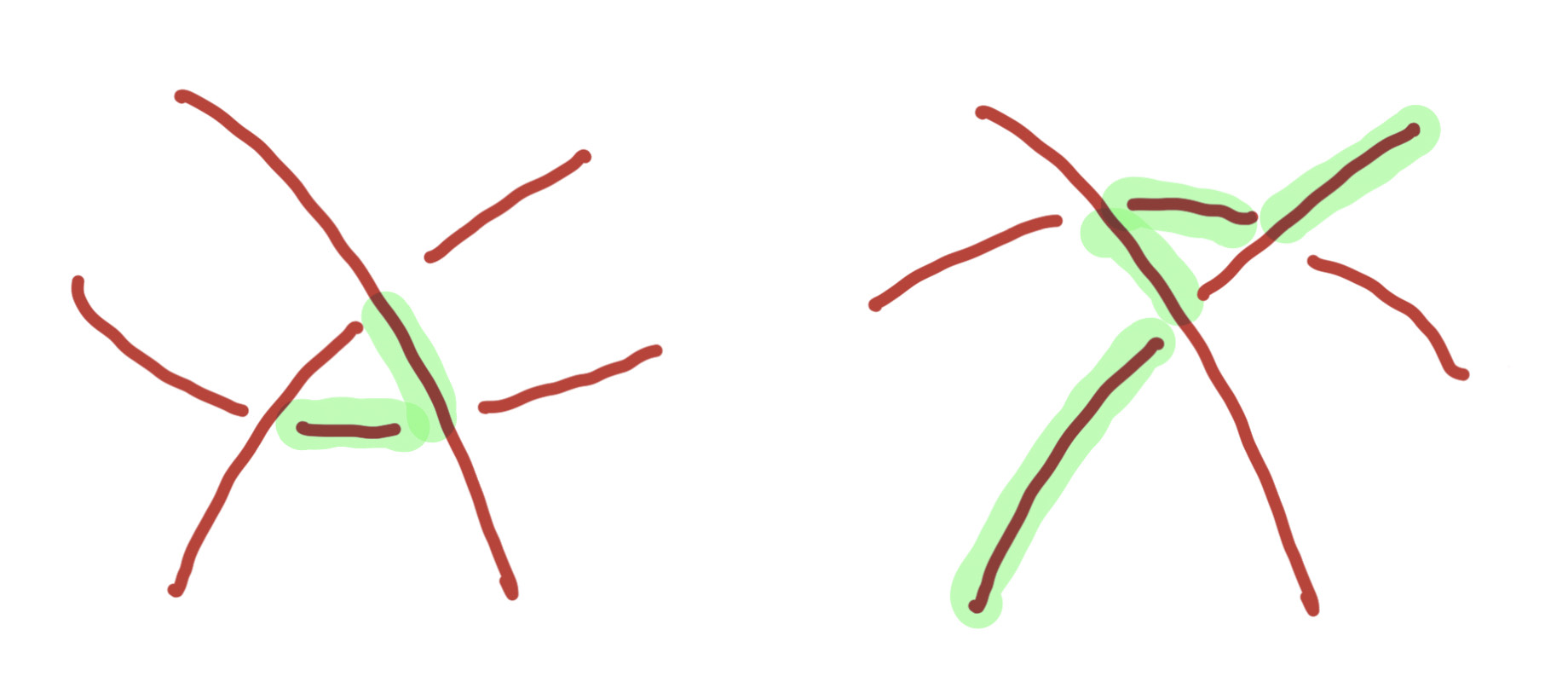}}
\put(150,70){$\Rightarrow$}
\end{picture}
\caption{The Reidemeister move 3 changes the number of non-alternating edges}\label{R3move}
\end{figure}

Figure~\ref{R3move} shows that the Reidemeister move 3 may increase the non-alternating edges by 2 without changing the number of crossings. 
We modified Dowker-Thistlethwaite codes of the remaining knots if the Reidemeister move 3 increases non-alternating edges, and then followed the same method as in the first stage, thus obtained 1,684 knots having arc index 13.

In the third stage, we've converted Dowker-Thistlethwaite codes of the remaining 2,990 knots without increasing the number of crossings. By the same method as above, we've obtained  2,467 knots with arc index 13.

In the last stage, for the remaining 523 knots, new tools were developed, which generate (some portion of) a list of grid diagrams (complete up to some size-preserving moves), and then track down therein grids matching a list of candidate knots.

\section{Minimal Grid Diagrams of Arc Index 13}

\begin{figure}[h!]
$$\begin{aligned}
\begin{tikzpicture}[scale=0.15]
\input tikz14n10ba.tex
\end{tikzpicture}
\begin{picture}(40,120)(10,10)
\put(20,70){$\Rightarrow$}
\end{picture}
&
\begin{tikzpicture}[scale=0.15]
\input tikz14n10ca.tex
\end{tikzpicture}\\
&\begin{tikzpicture}[scale=0.15]
\input tikz14n10c.tex
\end{tikzpicture}
\begin{picture}(40,120)(10,10)
\put(20,70){$\Rightarrow$}
\end{picture}
\begin{tikzpicture}[scale=0.15]
\input tikz14n10d.tex
\end{tikzpicture}
\begin{picture}(40,120)(10,10)
\put(20,70){$\Rightarrow$}
\end{picture}
\begin{tikzpicture}[scale=0.15]
\input tikz14n10e.tex
\end{tikzpicture}
\end{aligned}$$
\caption{Bottom 7 edges moved to the top, left 6 edges  moved to the right, then mirrored upside down}
\label{edge moves}
\end{figure}
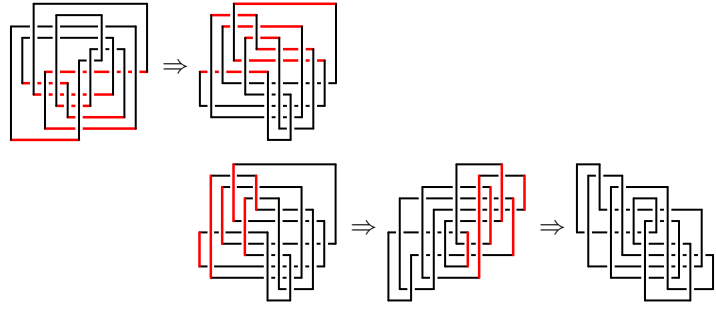

Grid diagrams have the property that outermost edges on one side can be shifted to the other side without changing the knot type, as shown in the first t and second steps of Figure~\ref{edge moves}. As the Dowker-Thistlethwaite code does not distinguish mirror images, we may reflect the diagram as shown in the third step of the figure.

After reducing the size from 16 to 13, we applied such edge moves and flips to make the diagrams to  look less random.
Due to a limit of pages, we only show the first 200 diagrams. 

\bigskip
\bigskip

{
\def\h#1#2#3{\draw[black,line width=0.75pt, cap=round] (#1,#2)--(#1+#3,#2);}
\def\w{0.2}
\def\v#1#2#3{
       \draw[white,line width=2.5pt] (#1,#2+\w)--(#1,#2+#3-\w);
       \draw[black,line width=0.75pt, cap=round] (#1,#2)--(#1,#2+#3);}
\def\vv#1#2#3{
       \draw[black,line width=0.75pt, cap=round] (#1,#2)--(#1,#2+#3);}
\def\knum#1#2#3#4{
\draw (-0.5,-0.5) (15.0,17);
\filldraw[black] (0,15.0) node[anchor=west]{$#1#2#3$}; }

\scriptsize

\noindent\input tikz8027_30_54.tex

}

\def\knum#1#2#3#4{
\draw (-0.5,-0.5) (15.0,17);
\filldraw[black] (0,15.7) node[anchor=west]{$#1#2#3$}; }

\section{The case of arc index 14}
We used the method of \emph{filtered spanning tree} as explained above, except that we did not need to change knot diagrams to increase the number of non-alternating edges from 4 to 6.
Due to a limit of pages, we only show the first 50 diagrams. 

\bigskip
{ \scriptsize \noindent
\input tikz15735red_36_54.tex }
\begin{tikzpicture}[scale=0.130]  
\hw11{12}
\end{tikzpicture} 
\begin{tikzpicture}[scale=0.130]  
\hw11{12}
\end{tikzpicture} 
\begin{tikzpicture}[scale=0.130]  
\hw11{12}
\end{tikzpicture} 
\begin{tikzpicture}[scale=0.130]  
\hw11{12}
\end{tikzpicture} 
\begin{tikzpicture}[scale=0.130]  
\hw11{12}
\end{tikzpicture}

{
\bigskip\color{blue}
\subsection*{Remark}
Due to a limit of pages, we only show the first 200 and 50  of the 8,027 and 17,535 diagrams, respectively. The full lists of minimal grid diagrams are available by typesetting the main source file with a little adjustment. Here are the steps.
\begin{enumerate}
\item Download all source files of this article from arXiv.
\item Remove \verb|`\endinput'| in the 5600th line of the file
 \verb#`tikz8027_30_54.tex'#
\item Remove \verb|`\endinput'| in the 1500th line of the file
 \verb#`tikz15735red_36_54.tex'#
\item Use \verb|pdfLaTeX| to typeset the  main  file 
\verb|`arc13-14xing14_arXiv.tex'|.
\item To remove this remark, remove \verb|`%'| 
 in the lines 559 and 576 of the main file before the step (4).
\end{enumerate}
}

\section*{Acknowledgments}
This work was supported by the Research \& Education Program at the   Korea Science Academy of KAIST with funds from the Ministry of Science and ICT.  
Lee and Stoimenow were supported by the National Research Foundation of Korea(NRF) grant funded by the Korea government(MSIT) (No. 2023R1A2C1003749).  
We used the program Knotscape at the beginning to get Dowker-Thistlethwaite codes of the 46,972 knots to program  for the filtered spanning trees and at the end to confirm the 8,027 and 15,735 grid diagrams represent the original ones~\cite{knotscape}. 
We also used the program KnotPlot to obtain the Dowker-Thistlethwaite codes from the final grid diagrams~\cite{knotplot}.
The knot diagrams of $14n10$ in Figure~\ref{fig:14n10edges} and $14n7951$ in Figure~\ref{R3move} were obtained from SnapPy~\cite{SnapPy}.

For the part of the arc index 13, the second group of authors did most of the programming to find spanning trees and for  the part of the arc index 14, the third group of authors did so. They worked independently in two different years.

\end{document}